\documentclass[10pt]{amsart}
\def\polhk{\c}
\font\msy=msbm10
\def\rtimes{\hbox{\msy o}}
\def\agm{{\rm agm}}
\def\sech{{\rm sech}}
\def\Sym{{\rm Sym}}
\def\Npow#1{{N^{(#1)}}}
\theoremstyle{plain}
\newtheorem{prop}{Proposition}[section]
\newtheorem{lemma}[prop]{Lemma}
\newtheorem{theorem}[prop]{Theorem}
\theoremstyle{remark}
\newtheorem{remark}[prop]{Remark}

\begin{document}
\title{Explicit lower bounds on the modular degree of an elliptic curve}
\author{Mark Watkins}
\address
{Department of Mathematics, McAllister Building,
The Pennsylvania State University, University Park, PA, 16802}

\begin{abstract}
We derive an explicit zero-free region for symmetric square \hbox{L-functions}
of elliptic curves, and use this to derive an explicit lower bound for the
modular degree of rational elliptic curves. The techniques are similar to
those used in the classical derivation of zero-free regions for Dirichlet
\hbox{L-functions}, but here, due to the work of Goldfield-Hoffstein-Lieman,
we know that there are no Siegel zeros, which leads to a strengthened result.
\end{abstract}
\maketitle

\section{Introduction}
Let $E$ be a rational elliptic curve of conductor~$N$;
by the work of Wiles~\cite{wiles} and others,
it is known that there is a surjective
map $\phi$ from $X_0(N)$ to~$E$ known as a modular parametrisation.
Our aim in this paper is to indicate sundry lower bounds on the modular
degree of an elliptic curve. Our starting point is a formula
of convolution type essentially
due to Shimura~\cite{shimura76} which states that we have
$${L(\Sym^2 E,1)\over 2\pi\Omega}=
{{\rm deg}\,\phi\over Nc^2}\prod_{p^2|N} U_p(1),$$
where $L(\Sym^2 E,s)$ is the motivic symmetric-square $L$-function of $E$
normalised so that $s=1/2$ is the point of symmetry, $\Omega$ is the area
of the fundamental parallelogram associated to the curve,
$c$ is the Manin constant which is known to be an integer
(see~\cite{edixhoven} or \cite[1.6]{stevens}),
and the $U_p(1)$ are fudge factors that can be given explicitly.
From the above we have that
$${\rm deg}\,\phi={Nc^2\over 2\pi\Omega}\cdot L(\Sym^2 E,1)\cdot
\prod_{p^2|N} U_p(1)^{-1}.$$
One of our goals is to show a bound
of the type ${\rm deg}\,\phi\gg N^{7/6-\epsilon}$ as $N\rightarrow\infty$.
Indeed, this has been known in folklore (see for instance the paper
by Papikian~\cite{papikian} that deals with a functional field analogue)
since the time of the Goldfeld-Hoffstein-Lieman \cite{GHL} appendix
to the work of Hoffstein-Lockhart~\cite{HL},
but herein we give a more complete
proof and compute explicit constants. We shall assume that $N\ge 20000$,
(and thus we have that the symmetric-square conductor $\Npow2\ge 142$),
as else the tables of Cremona~\cite{cremona} can be used.

Previously, explicit bounds had been obtained in a couple of ways.
As N.~Elkies pointed out to us, one can use an idea of Ogg~\cite{ogg}
to show that $d={\rm deg}\,\phi \gg N/p$ where $p$ is any prime of
good reduction. Here is the argument.
Reduce the modular parametrisation map mod $p$,
and consider it over the field $k$ of $p^2$ elements.
Now $X_0(N)$ has about $pN/12$ supersingular points, all defined over~$k$;
whereas the elliptic curve has at most $(p+1)^2$ $k$-rational points.
Since each of these has at most $d$ preimages, and these preimages
must include all the $k$-rational points, the estimate  $d\gg N/p$ follows,
with a constant of $1/12$ as $N\rightarrow\infty$.

One can make a ``characteristic zero'' version of this argument
by using lower bounds for the eigenvalues of the Laplacian on $X_0(N)$
in place of the supersingular points. This allows one to obtain a linear
lower bound on the modular degree.
Already in a paper of Li and Yau~\cite{liyau}
there appears the technique for passing from an eigenvalue bound
to a modular degree bound, and this is then made explicit
by Abramovich~\cite{abramovich}
whose result is equivalent to ${\rm deg}\,\phi \ge 7N/1600$
which would improve to ${\rm deg}\,\phi \ge N/192$,
upon assuming the Selberg eigenvalue conjecture.
We shall see that these linear bounds are superior to ours
unless $N$ is quite large, as the constant in our bound
will turn out to be rather small.

\section{Bounding the area of the fundamental parallelogram}
Let $E$ be an elliptic curve, which we write in the form
$y^2=4x^3+b_2x^2+2b_4x+b_6$ in such a way that the discriminant
is minimal away from~2. The polynomial on the right side of this equation
is sometimes called the 2-torsion polynomial.
There are two natural cases depending on whether
the discriminant is positive or negative.
In either case we have that $\Omega$ is the real period multiplied by
the imaginary part of the imaginary period.
We have the following lower bound on $1/\Omega$.

\begin{lemma}
Let $E$ be an elliptic curve, $\Omega$ the area of its fundamental
parallelogram, and $D$ the absolute value of its discriminant.
Then $1/\Omega \ge {D^{1/6}\over 14.045}$.
\end{lemma}

\begin{proof}
The proof naturally divides into 2 cases, depending on whether $\Delta>0$.

{\bf Case I: positive discriminant.}
When the discriminant is positive
the 2-torsion polynomial has three real roots,
which we order as $e_1>e_2>e_3$.
We then have that (see Chapter 7 of Cohen~\cite{cohen})
the real period of $E$ is $\pi/\agm(\sqrt{e_1-e_2},\sqrt{e_1-e_3})$
and the imaginary period is
$\pi i/\agm(\sqrt{e_2-e_3},\sqrt{e_1-e_3})$,
and that we also have $\sqrt{\Delta/16}=(e_1-e_2)(e_1-e_3)(e_2-e_3)$.
Let $t={e_1-e_2\over e_1-e_3}$ so that we have $t\in (0,1)$
and $(e_1-e_3)\cdot [4t(1-t)]^{1/3}=\Delta^{1/6}$,
and recall that $\agm(x,y)=x\cdot\agm(1,y/x)$, implying
\begin{align*}
1/\Omega&={1\over\pi^2}(e_1-e_3)\cdot
\agm\bigl(1,\sqrt t\bigr)\cdot\agm\bigl(1,\sqrt{1-t}\bigr)\\
&\ge {1\over \pi^2}(e_1-e_3)
\cdot[4t(1-t)]^{1/3}\cdot\agm\bigl(1,1/\sqrt 2\bigr)^2
={D^{1/6}\over\pi^2}\cdot\agm\bigl(1,1/\sqrt 2\bigr)^2,
\end{align*}
where the inequality follows from calculus,
the relevant quotient function being minimised at $t=1/2$.

{\bf Case II: negative discriminant.}
When the discriminant is negative we let $r$ be the real root of
the 2-torsion polynomial,
and write $\tilde r=r+b_2/12$, so that $-\tilde r/2\pm iZ$
are the other roots. The real period is now
$2\pi/\agm(2\sqrt B,\sqrt{2B+A})$
and the (vertical part of the) imaginary period is
$\pi i/\agm(2\sqrt B,\sqrt{2B-A})$ where $A=3r+b_2/4=3\tilde r$ and
$B=\sqrt{3r^2+b_2r/2+b_4/2}=\sqrt{(3\tilde r/2)^2+Z^2}$.
Also note that $2ZB^2=\sqrt{-\Delta/16}$.
Write $c=\tilde r/Z$, so that $A=3cZ$ and $B=Z\sqrt{1+9c^2/4}$,
so that $D^{1/6}=2Z(1+9c^2/4)^{1/3}$. Writing $M(x)=\agm(1,x)$, we have
\begin{align*}
1/\Omega&={1\over 2\pi^2}\cdot
2\sqrt B\cdot\agm\biggl(1,\sqrt{2B+A\over 4B}\biggr)\cdot
2\sqrt B\cdot\agm\biggl(1,\sqrt{2B-A\over 4B}\biggr)\\
&={1\over 2\pi^2}\cdot 4Z\sqrt{1+9c^2/4}\cdot
M\Biggl(\sqrt{{1\over 2}+{3c\over \sqrt{16+36c^2}}}\Biggr)\cdot
M\Biggl(\sqrt{{1\over 2}-{3c\over \sqrt{16+36c^2}}}\Biggr)\\
&={D^{1/6}\over\pi^2}\cdot(1+9c^2/4)^{1/6}\cdot
M\Biggl(\sqrt{{1\over 2}+{3c\over \sqrt{16+36c^2}}}\Biggr)\cdot
M\Biggl(\sqrt{{1\over 2}-{3c\over \sqrt{16+36c^2}}}\Biggr)\\
&\ge{D^{1/6}\over\pi^2}\cdot4^{1/6}\cdot
\agm\Biggl(1,\sqrt{{1\over 2}+{\sqrt 3\over 4}}\Biggr)\cdot
\agm\Biggl(1,\sqrt{{1\over 2}-{\sqrt 3\over 4}}\Biggr),
\end{align*}
as the function is minimised at $c=\pm\sqrt{4/3}$.
In both cases we have $1/\Omega \ge {D^{1/6}\over 14.045}$.
\end{proof}

\section{Zero-free regions and lower bounds for symmetric square L-functions}
We next turn to making the argument of \cite{GHL} explicit.
We first need to derive a zero-free region for $L(\Sym^2 E,s)$,
and then turn this into a lower bound for $L(\Sym^2 E,1)$.
\subsection{Zero-free regions for curves without complex multiplication}
\begin{lemma}
Let $L(\Sym^2 f_E,s)$ be the symmetric-square $L$-function of $f_E$,
where $f_E$ is the form associated to a rational elliptic curve~$E$
that does not have complex multiplication.
Then $L(\Sym^2 f_E,s)$ has no real zeros with $s\ge 1-\delta/\log(\Npow2/C)$,
where $\delta=2(5-2\sqrt 6)/5\approx 0.040408$, $C=96$, and $\Npow2\ge 142$
is the symmetric-square conductor of~$E$.
\end{lemma}

\begin{proof}
We follow the proof in the appendix \cite{GHL} of~\cite{HL},
which uses the idea that a function with a double pole at $s=1$
cannot have a triple zero too close to this pole.
The product \hbox{$L$-function} in question (see~page~180) is
${\bf L}(s)=\zeta(s)\cdot L(F,s)^3\cdot L(\Sym^2 F,s)
=\zeta(s)^2\cdot L(\Sym^2 f_E,s)^3\cdot L(\Sym^4 f_E,s)$ where $F=\Sym^2 f_E$
and all \hbox{$L$-functions} of symmetric powers are motivic.
As that paper notes earlier in a slightly different context
(see page 167, after the proof of Lemma~1.2), we have that
the Dirichlet series ${\bf L}(s)$ has nonnegative coefficients at
primes of good reduction, and, more important for our immediate purposes,
by taking the logarithmic derivative we see that $({\bf L}^\prime/{\bf L})(s)$
has nonpositive coefficients at such primes.
It is asserted in \cite{HL} that the Langlands correspondence implies the
nonpositivity at bad primes.
For our case of elliptic curves, the proper Euler factor at bad primes
is worked out in the Sheffield dissertation of Phil Martin \cite{martin}, and
it can be verified directly that we do indeed have the desired nonpositivity.
Note that D\polhk abrowski \cite{dabrowski}
claims to compute the Euler factors at bad primes
in Lemma 1.2.3 on page 63 of that paper,
but the method used therein appears to be erroneous;
similarly the method in an appendix of a previous version of this paper
failed to consider the cases of noncyclic inertia group correctly.

We also need to compute the factor at infinity
and the conductor of ${\bf L}(s)$.
For the factor at infinity, this is done on pages~60--61 of~\cite{dabrowski}:
we have a factor of $\Gamma(s/2)/\pi^{s/2}$ for $\zeta(s)$,
a factor of $\Gamma(s+1)\Gamma\bigl((s+1)/2\bigr)/(4\pi^3)^{s/2}$
for $L(\Sym^2 f_E,s)$, and a factor of
$\Gamma(s+2)\Gamma(s+1)\Gamma(s/2)/(16\pi^5)^{s/2}$ for $L(\Sym^4 f_E,s)$.
For the bad primes, we refer to \cite{martin};
we should note that we can bound the symmetric-square conductor
by the square of the conductor, that is, $\Npow2\le N^2$,
and similarly the symmetric-fourth-power conductor is bounded
by the square of the symmetric-square conductor, that is,
$\Npow4\le (\Npow2)^2$. This also follows from~\cite{BH}.
Note the symmetric-square conductor is actually a square,
and so some authors (for instance~\cite{watkins})
define it to be the square root of our choice here.

So we claim that
\begin{eqnarray*}
\Phi(s)&=\Gamma(s/2)^3\Gamma(s+1)^4\Gamma\bigl((s+1)/2\bigr)^3\Gamma(s+2)
\Biggl(\displaystyle{{\Npow2}^3\Npow4\over 1024\pi^{16}}\Biggr)^{s/2}\cdot\\
&\zeta(s)^2\cdot L(\Sym^2 f_E,s)^3\cdot L(\Sym^4 f_E,s)
\end{eqnarray*}
is meromorphic and symmetric under the map $s\rightarrow 1-s$.
The asserted analytic properties follow from work of
Gelbart and Jacquet \cite{GJ}
and Shimura \cite{shimura75} for the symmetric square, and later authors such
as Kim and Shahidi for higher symmetric powers~\cite{KS}.
By Bump and Ginzburg~\cite{BG}, when $f_E$ is not a $GL(1)$-lift
(when $E$ does not have complex multiplication),
$\Phi(s)$ has a double pole at $s=1$ (see also the work of Kim~\cite{kim}).

So the function $\Lambda(s)=s^2(1-s)^2\Phi(s)$ is entire,
and by taking the logarithmic derivative of its Hadamard product,
we get that
$\sum_{\rho} {w_\rho\over s-\rho}=
{2\over s-1}+{2\over s}+{\Phi^\prime\over\Phi}(s)$,
where $w_\rho$ is an appropriate weight for the zero
and the order of summation is taken over
conjugate pairs of zeros of $\Phi(s)$, and so is convergent.
\smallskip
Now assume that $L(\Sym^2 f_E,s)$ has a zero at $\beta$.
Then $\Phi(s)$ has a triple zero at $\beta$, so that we have
\begin{align*}
{3\over s-\beta}+&{3\over s-(1-\beta)}+\sum_{\rho} {w_\rho\over s-\rho}=\\
&{2\over s-1}+{2\over s}+
{3\over 2}\log\Npow2+{1\over 2}\log\Npow4-\log 32\pi^8+\\
&\kern48pt+3{\Gamma^\prime\over\Gamma}(s/2)+4{\Gamma^\prime\over\Gamma}(s+1)+
3{\Gamma^\prime\over\Gamma}((s+1)/2)+{\Gamma^\prime\over\Gamma}(s+2)
+{{\bf L}^\prime\over{\bf L}}(s),
\end{align*}
where the sum over $\rho$ is over the non-Siegel zeros of $\Phi(s)$.

Now assume that $L(\Sym^2 f_E,s)$ has a zero at
$\beta\ge 1-2(5-2\sqrt 6)/5\log(\Npow2/C)$.
We let $C=96$ and write $\delta=(1-\beta)\log(\Npow2/C)$ and
evaluate the above displayed equation at
$s=\sigma=1+\eta\delta/\log (\Npow2/C)$ where
$\eta={1\over 10\delta}\bigl[(2-5\delta)-\sqrt{25\delta^2-100\delta+4}\bigr]$
is the smaller positive root of ${5\over 2}\delta x^2+({5\over 2}\delta-1)x+2$.
Note that both roots are real and positive when $0<\delta\le 2(5-2\sqrt 6)/5$.
We get a crude lower bound of zero for the $\rho$-sum
by pairing conjugate roots, and so
\begin{align*}
{3\over \sigma-\beta}\le
{2\over \sigma-1}+&{2\over \sigma}-{3\over \sigma-(1-\beta)}+
{5\over 2}\log\Npow2-\log 32\pi^8+\\
&+3{\Gamma^\prime\over\Gamma}(\sigma/2)+
4{\Gamma^\prime\over\Gamma}(\sigma+1)+
3{\Gamma^\prime\over\Gamma}((\sigma+1)/2)+
{\Gamma^\prime\over\Gamma}(\sigma+2).
\end{align*}
From this we get that
\begin{align*}
{3\over\eta+1}&{\log (\Npow2/C)\over\delta}\le\\
&{2\log (\Npow2/C)\over\eta\delta}+
{2\over \sigma}-{3\over \sigma-(1-\beta)}+
{5\over 2}\log(\Npow2/C)-\log 32\pi^8+\\
&+3{\Gamma^\prime\over\Gamma}(\sigma/2)
+4{\Gamma^\prime\over\Gamma}(\sigma+1)
+3{\Gamma^\prime\over\Gamma}((\sigma+1)/2)
+{\Gamma^\prime\over\Gamma}(\sigma+2)+{5\over 2}\log C,
\end{align*}
and here the terms with $\log(\Npow2/C)$ cancel due to the definition
of~$\eta$. So we have
\begin{align*}
0\le {2\over \sigma}-&{3\over \sigma-(1-\beta)}-\log 32\pi^8+\\
&+3{\Gamma^\prime\over\Gamma}(\sigma/2)+
4{\Gamma^\prime\over\Gamma}(\sigma+1)+
3{\Gamma^\prime\over\Gamma}((\sigma+1)/2)+
{\Gamma^\prime\over\Gamma}(\sigma+2)+{5\over 2}\log C,
\end{align*}
Now $\eta\delta$ is maximised at the endpoint where $\delta=(5-2\sqrt 6)/5$,
giving us that $\sigma\le 1+2(\sqrt 6-2)/5\log(\Npow2/C)$.
Under our assumption that $\Npow2\ge 142$ and definition of $C=96$,
this gives that $\sigma\le 1.46$, so that
$$3{\Gamma^\prime\over\Gamma}(\sigma/2)+4{\Gamma^\prime\over\Gamma}(\sigma+1)+
3{\Gamma^\prime\over\Gamma}((\sigma+1)/2)+{\Gamma^\prime\over\Gamma}(\sigma+2)
\le 1.74.$$
We also have that
\begin{align*}
{2\over \sigma}&-{3\over \sigma-1+\beta}=
{2\over 1+\eta\delta/\log(\Npow2/C)}-
{3\over 1+\delta(\eta-1)/\log(\Npow2/C)}\\
&\le {2\over 1+2(\sqrt 6-2)/5\log(\Npow2/C)}-
{3\over 1+2(-7+3\sqrt 6)/5\log(\Npow2/C)}\le -0.84,
\end{align*}
so that we get the contradiction that
$0\le -0.84-12.62+1.74+{5\over 2}\log C\le -0.30$.
Thus there are no zeros in the region indicated.
\end{proof}

\begin{remark}
The constant $\delta$ can be improved if we could lower-bound
$\sum_\rho {1\over s-\rho}$ less crudely as some constant times
$\log(\Npow2)$, which is likely feasible by zero-counting arguments.
The constant $C$ can be improved simply by requiring $\Npow2$ to be larger.
\end{remark}

\subsection {Zero-free regions for curves without complex multiplication}
\begin{lemma}
Let $E$ be a rational elliptic curve with complex multiplication by an order
in the complex quadratic field~$K$. Then $L(\Sym^2 f_E,s)$ has no real zeros
with $\sigma\ge 1-\delta/\log(\Npow2/C)$,
where here we have $\delta=2^{1/2}+2-2^{7/4}\approx 0.050628$,
$C=64$, and $\Npow2\ge 142$
is the symmetric-square conductor of~$E$.
\end{lemma}

\begin{proof}
When $E$ has complex multiplication by an order of~$K$,
the representation associated to $f_E$ is dihedral,
and so by \cite{kim} the fourth symmetric power $L$-function
has a pole at~$s=1$, so that the $\Phi(s)$ of above has a triple pole at $s=1$.
However, as noted by~\cite{GHL}, in this case we have that $L(\Sym^2 f_E,s)$
can be factored.
Recall that there is some Hecke character $\psi$ of~$K$ such that
$L(f_E,s)=L(\psi,K,s)$ with $\psi(z)=\chi(z)(z/|z|)$ for some
character~$\chi$ defined on the ring of integers of~$K$.
Here $\chi$ has order at most 6, and is of order 1 or 2
unless $K$ is ${\bf Q}(i)$ or ${\bf Q}(\zeta_3)$.
We have the factorisation
$L(\Sym^2 f_E,s)=L(\theta_K,s)L(\psi^2,K,s)$ where $\theta_K$ is the
quadratic character of the imaginary quadratic field~$K$.
Here $\psi^2$ is the ``motivic'' square of $\psi$,
so that if $\psi(z)=\chi(z)(z/|z|)$ for some quadratic character~$\chi$,
we then have $\psi^2(z)=(z/|z|)^2$. Thus the square of $\chi$ is the
trivial character on $K$ and not the principal character of the same
modulus of $\chi$. The same convention shall apply to higher symmetric powers.

Similar to the above factorisation of the symmetric-square $L$-function,
by comparison of Euler factors we find that
$L(\Sym^4 f_E,s)=\zeta(s)L(\psi^2,K,s)L(\psi^4,K,s)$ and
$L(\Sym^6 f_E,s)=L(\theta_K,s)L(\psi^2,K,s)L(\psi^4,K,s)L(\psi^6,K,s)$
Here we can note that $L(\Sym^4 f_E,s)$ has a pole at $s=1$
but $L(\Sym^6 f_E,s)$ does not.

For the seven choices of $K$ with
${\rm disc}(K)<-4$, we thus have only one function $L(\psi^2,K,s)$ to consider,
and a direct computation establishes the indicated zero-free region.
For $K={\bf Q}(i)$ we need to consider quartic twists, and for
$K={\bf Q}(\zeta_3)$ we need to consider both cubic and sextic twists.
Note that Theorem~2 of Murty \cite{murty} erroneously only considers
quadratic twists, and thus the proof that the modular degree is at
least $N^{3/2-\epsilon}$ for elliptic curves with complex multiplication
is wrong. In fact, simply by taking sextic twists of $X_0(27)$
we can easily achieve a growth rate of only $N^{7/6+\epsilon}$.

{\bf Case I:}
We first consider the case where $K={\bf Q}(i)$.
Using the above decomposition of the symmetric-square $L$-function, we get
that the completed $L$-function that is symmetric under $s\rightarrow 1-s$ is
$$\biggl({\Npow2/4\over 4\pi^2}\biggr)^{s/2}\Gamma(s+1)L(\psi^2,K,s).$$
In order for the fourth symmetric power to work out, we see that
$$\biggl({4\Npow4/\Npow2\over 4\pi^2}\biggr)^{s/2}\Gamma(s+2)L(\psi^4,K,s)$$
is symmetric under $s\rightarrow 1-s$. Here we have $\Npow4=\Npow2$
from~\cite{martin}, due to the fact that the relevant
inertia groups are all $C_2$, $C_4$, or~$Q_8$.

The standard ingredient of proofs of a zero-free region for a
Hecke $L$-function is a trigonometric polynomial that is always nonnegative.
Here we take $(1+\sqrt 2\cos\theta)^2=2+2\sqrt 2\cos\theta+\cos 2\theta$.
A better result might come about from using higher degree cosine
polynomials, but the $\Gamma$-factors might be burdensome.
Also note that the work of Coleman \cite{coleman} could be used
if we did not need to be explicit.
Note that at bad primes we still have the desired nonpositivity
since $2\ge\cos 2\theta$ for all~$\theta$.

So we are led to consider the nonpositive sum
$${{\bf L}^\prime\over{\bf L}}(s)
=2{\zeta^\prime\over\zeta}(s)+2\sqrt2{L^\prime\over L}(\psi^2,K,s)+
{L^\prime\over L}(\psi^4,K,s).$$
Assume there is a zero of $L(\psi^2,K,s)$ at~$\beta$.
By the functional equation we get
\begin{align*}
{2\sqrt 2\over s-\beta}+&{2\sqrt 2\over s-(1-\beta)}+
\sum_{\rho} {w_\rho\over s-\rho}=\\
&{2\over s-1}+{2\over s}+
2\log(1/\sqrt\pi)+2\sqrt 2\log(\sqrt{\Npow2}/4\pi)+\log(2/2\pi)+\\
&\kern72pt
+2{\Gamma^\prime\over\Gamma}(s/2)+2\sqrt 2{\Gamma^\prime\over\Gamma}(s+1)+
{\Gamma^\prime\over\Gamma}(s+2)+{{\bf L}^\prime\over {\bf L}}(s),
\end{align*}
where, in the sum over zeros, $w_\rho$ is an appropriate weight for the zero.

Now we assume that the function $L(\Sym^2 f_E,s)=L(\theta_K,s)L(\psi^2,K,s)$
has a zero at $\beta\ge 1-(2^{1/2}+2-2^{7/4})/\log(\Npow2/C)$.
We define $C=100$ and write $\delta=(1-\beta)\log(\Npow2/C)$ and
proceed to evaluate the above displayed equation at
$s=\sigma=1+\eta\delta/\log (\Npow2/C)$ where $\eta$ is given by
the smaller positive root of $\delta\sqrt 2 x^2+(\delta\sqrt 2-2\sqrt 2+2)x+2$.
Note that both roots are real and positive
when $0<\delta\le 2^{1/2}+2-2^{7/4}$.

We again get a crude lower bound of zero for the $\rho$-sum
by pairing conjugate roots and have that
$({\bf L}^\prime/{\bf L})(\sigma)\le 0$, and so
\begin{align*}
{2\sqrt2\over \sigma-\beta}\le&
{2\over \sigma-1}+{2\over \sigma}-{2\sqrt 2\over \sigma-(1-\beta)}
+\log(1/\pi)+\sqrt 2\log\Npow2+\\
&\kern6pt
+2\sqrt 2\log(1/4\pi)+\log(1/\pi)+
+2{\Gamma^\prime\over\Gamma}(\sigma/2)+
2\sqrt2{\Gamma^\prime\over\Gamma}(\sigma+1)+
{\Gamma^\prime\over\Gamma}(\sigma+2).
\end{align*}

From this we get that
\begin{align*}
{2\sqrt 2\over\eta+1}{\log (\Npow2/C)\over\delta}
\le{2\log (\Npow2/C)\over\eta\delta}&+
{2\over \sigma}-{2\sqrt 2\over \sigma-(1-\beta)}+\sqrt 2\log(\Npow2/C)+\\
&+2\log(1/\pi)+2\sqrt 2\log(1/4\pi)+
+2{\Gamma^\prime\over\Gamma}(\sigma/2)+\\
&+2\sqrt2{\Gamma^\prime\over\Gamma}(\sigma+1)+
{\Gamma^\prime\over\Gamma}(\sigma+2)+\sqrt 2\log C
\end{align*}
and here the terms with $\log(\Npow2/C)$ cancel due to the definition
of~$\eta$. So we have
\begin{align*}
0\le {2\over \sigma}-{2\sqrt 2\over \sigma-(1-\beta)}
&+2\log(1/\pi)+2\sqrt 2\log(1/4\pi)+\\
&+2{\Gamma^\prime\over\Gamma}(\sigma/2)+
2\sqrt2{\Gamma^\prime\over\Gamma}(\sigma+1)+
{\Gamma^\prime\over\Gamma}(\sigma+2)+\sqrt 2\log C.
\end{align*}
Now $\delta\eta$ is maximised as $\sqrt 2(2^{1/4}-1)$
when $\delta=2^{1/2}+2-2^{7/4}$, and so under our assumption that
$\Npow2\ge 142$ and $C=100$ we have that $\sigma\le 1.8$,
so that the $\Gamma$-terms contribute less than 2.821.
We also have that
\begin{align*}
{2\over \sigma}-{2\sqrt 2\over \sigma-1+\beta}&=
{2\over 1+\eta\delta/\log(\Npow2/C)}-
{2\sqrt 2\over 1+\delta(\eta-1)/\log(\Npow2/C)}\\
&\kern-32pt\le {2\over 1+\sqrt 2(2^{1/4}-1)/\log(\Npow2/C)}-
{2\sqrt 2\over 1+(3\cdot 2^{3/4}-2\sqrt 2-2)/\log(\Npow2/C)}\\
&\kern-32pt\le -0.612,
\end{align*}
so that we get the contradiction that
$0\le -0.612-9.448+2.821+\sqrt 2\log C\le -0.726$.

{\bf Case II:}
Next we consider the other case where $K={\bf Q}(\zeta_3)$.
By examining the above functional equations for
symmetric-power $L$-functions, we find that
$$\biggl({\Npow2/3\over 4\pi^2}\biggr)^{s/2}\Gamma(s+1)L(\psi^2,K,s),
\>\>\biggl({3\Npow4/\Npow2\over 4\pi^2}\biggr)^{s/2}\Gamma(s+2)L(\psi^4,K,s),$$
$${\rm and}\>\>\>
\biggl({\Npow6/3\Npow4\over 4\pi^2}\biggr)^{s/2}\Gamma(s+3)L(\psi^6,K,s)$$
are all invariant under $s\rightarrow 1-s$.
Using~\cite{martin}, the fact that all the relevant inertia groups
are $C_3$ or $C_6$ implies that we have that $\Npow6=\Npow4=(\Npow2)^2$
except in the case when $3^3\|N$, when the inertia group is the
semi-direct product $C_3\rtimes C_4$ and we have
$\Npow6=9\Npow4=(\Npow2)^2$.

Here we choose a trigonometric polynomial of the form
$(1+\cos\theta)(1+\beta\cos\theta)^2$. It turns out that
the optimal $\beta$ for our purposes is twice the positive
root of the polynomial $x^5-25x^4-4x^3+30x^2+19x+3$, approximately 2.629152166,
but we do not lose much by taking $\beta=5/2$, so that our
nonnegative trigonometric polynomial is
${1\over 16}(106+171\cos\theta+90\cos 2\theta+25\cos 3\theta)$.
So we are led to consider the nonpositive sum
$${{\bf L}^\prime\over{\bf L}}(s)
=106{\zeta^\prime\over\zeta}(s)+171{L^\prime\over L}(\psi^2,K,s)+
90{L^\prime\over L}(\psi^4,K,s)+25{L^\prime\over L}(\psi^6,K,s),$$
with the nonpositivity at bad primes following as before.

Assume there is a zero of $L(\psi^2,K,s)$ at~$\beta$.
By the functional equation we get
\begin{align*}
{171\over s-\beta}&+{171\over s-(1-\beta)}+
\sum_{\rho} {w_\rho\over s-\rho}=\\
&\kern-6pt
{106\over s-1}+{106\over s}+106\log(1/\sqrt\pi)+
{171\over 2}\log(\Npow2/12\pi^2)+\\
&\kern-6pt
+45\log(\Npow4/\Npow2)+45\log(3/4\pi^2)
+{25\over 2}\log(\Npow6/\Npow4)+{25\over 2}\log(1/12\pi^2)+\\
&\kern-6pt
+106{\Gamma^\prime\over\Gamma}(s/2)+171{\Gamma^\prime\over\Gamma}(s+1)+
90{\Gamma^\prime\over\Gamma}(s+2)+25{\Gamma^\prime\over\Gamma}(s+3)+
{{\bf L}^\prime\over {\bf L}}(s).
\end{align*}

Now assume that $\beta\ge 1-(554-12\sqrt{2014})/261\log(\Npow2/C)$.
We let $C=64$ and $\delta=(1-\beta)\log(\Npow2/C)$ and
proceed to evaluate the above displayed equation at
$s=\sigma=1+\eta\delta/\log (\Npow2/C)$ where $\eta$ is given
by the smaller positive root of $261\delta x^2+(261\delta-130)x+212 $.
Note that both roots are real and positive
when $0<\delta\le (554-12\sqrt{2014})/261$.

We again get a crude lower bound of zero for the $\rho$-sum
by pairing conjugate roots and have that
$({\bf L}^\prime/{\bf L})(\sigma)\le 0$, and so
\begin{align*}
{171\over \sigma-\beta}\le
{106\over \sigma-1}&+{106\over \sigma}-{171\over \sigma-(1-\beta)}
+{261\over 2}\log\Npow2+\\
&+53\log(1/\pi)-171\log\sqrt{12}+171\log(1/\pi)+\\
&+45\log 3/4+90\log(1/\pi)-25\log\sqrt{12}+25\log(1/\pi)+\\
&+106{\Gamma^\prime\over\Gamma}(s/2)+171{\Gamma^\prime\over\Gamma}(s+1)+
90{\Gamma^\prime\over\Gamma}(s+2)+25{\Gamma^\prime\over\Gamma}(s+3)
\end{align*}

From this and the definition of $\eta$ we get
\begin{align*}
0\le {106\over \sigma}-{171\over \sigma-(1-\beta)}&
+339\log(1/\pi)-\log (3^{53}2^{286})+{261\over 2}\log C\\
&+106{\Gamma^\prime\over\Gamma}(s/2)+171{\Gamma^\prime\over\Gamma}(s+1)+
90{\Gamma^\prime\over\Gamma}(s+2)+25{\Gamma^\prime\over\Gamma}(s+3)
\end{align*}
Now $\delta\eta$ is maximised as $(6\sqrt{2014}-212)/261$
when $\delta=(554-12\sqrt{2014})/261$, and so under our assumption that
$\Npow2\ge 142$ and $C=64$ we have that $\sigma\le 1.28$,
so that the above $\Gamma$-sum is less than~153.
We thus get the contradiction that $0\le -59-644+543+153$.
\end{proof} 

\subsection{Lower bounds from zero-free regions}
We use these zero-free regions to lower bound~$L(\Sym^2 f_E,1)$.

\begin{lemma}
Let $E$ be a rational elliptic curve with whose symmetric-square
conductor satisfies $\Npow2\ge 142$.
Then $L(\Sym^2 f_E,1)\ge {0.033\over\log \Npow2}$.
\end{lemma}

\begin{proof}
We use Rademacher's formulation \cite{rademacher}
of the Phragm\'en--Lindel\"of Theorem
to bound $L(\Sym^2 f_E,s)$ and $\zeta(s)$.
First note that by the Euler product we have
$|L(\Sym^2 f_E,3/2+it)|\le \zeta(3/2)^3$,
and so by the functional equation we get
\begin{align*}
|L(\Sym^2 f_E,&-1/2+it)|=\\
&=(\Npow2/4\pi^3)
{|\Gamma(5/2+it)|\cdot|\Gamma(5/4+it/2)|\over
|\Gamma(1/2+it)|\cdot|\Gamma(1/4+it/2)|}
|L(\Sym^2 f_E,3/2+it)|\\
&\le\zeta(3/2)^3{\Npow2\over 4\pi^3}\cdot
\biggl|{3\over 2}+it\biggr|\cdot\biggl|{1\over 2}+it\biggr|\cdot
\biggl|{1\over 4}+{it\over 2}\biggr|
\le\zeta(3/2)^3{\Npow2\over 8\pi^3}\cdot
\biggl|{3\over 2}+it\biggr|^3.
\end{align*}
So using the result of Rademacher with
$Q=2$, $a=-1/2$, $b=3/2$, $\alpha=3$, $\beta=0$,
$A=\zeta(3/2)^3(\Npow2/8\pi^3)$, and $B=\zeta(3/2)^3$,
we get that
$$|L(\Sym^2 f_E,1/2+it)|\le
(A|5/2+it|^3)^{1/2} B^{1/2}
=\zeta(3/2)^3\cdot\sqrt{\Npow2\over8\pi^3}|5/2+it|^{3/2},$$
which is not anywhere near an optimal bound, but will suffice.
Similarly we get that
$|\zeta(1/2+it)|\le {\zeta(3/2)\over\sqrt{2\pi}}\sqrt{{9\over 4}+t^2}$.

Let $b=1-{1\over 25\log\Npow2}$ so that
${\bf L}(s)=L(\Sym^2 f_E,s)\zeta(s)$ has no zeros in $[b,1)$,
so that ${\bf L}(b)<0$. Note also that $b\ge 0.99$ due to our assumption
that $\Npow2\ge 142$. Writing ${\bf L}(s)=\sum_n a_n/n^s$ as a Dirichlet
series with nonnegative coefficients, by the Mellin transform we have that
$$\sum_n {a_n\over n^b} e^{-n/X}=
\int_{(2)} \Gamma(s) X^s {\bf L}(s+b) {ds\over 2\pi i}.$$
Via moving the contour to where the real part of $s$ is $1/2-b$,
we get that the integral is $RX^{1-b}\Gamma(1-b)+{\bf L}(\beta)+E(X)$
where $R=L(\Sym^2 f_E,1)$ and the error term $E(X)$ is bounded by
${1\over\pi}\int_0^\infty
\bigl|\Gamma(1/2-b+it)X^{1/2-b}{\bf L}(1/2+it)\bigr| dt$.
By another theorem in Rademacher \cite{rademacher} we have that
\begin{align*}
|\Gamma(1/2-b+it)|&\le |1/2+it|^{1-b}\cdot |\Gamma(-1/2+it)|\\
&\le {|1+it|^{0.01}\over |1/2+it|} |\Gamma(1/2+it)|
={2(1+t^2)^{1/200}\over \sqrt{1+4t^2}}\sqrt{\pi\sech \pi t}
\end{align*}
We compute that
$${\zeta(3/2)^4\over 4\pi^2}\int_0^\infty
\biggl({25\over 4}+t^2\biggr)^{3/2}\sqrt{{9\over 4}+t^2}\cdot
{2(1+t^2)^{1/200}\over \sqrt{1+4t^2}}\sqrt{\pi\sech \pi t}\, dt<62,$$
so that $|E(X)|\le 20\sqrt{\Npow2}\cdot X^{1/2-b}$.
Since ${\bf L}(\beta)\le 0$ we ge
$$\sum_n {a_n\over n^b} e^{-n/X}\le
RX^{1-b}\Gamma(1-b)+20\sqrt{\Npow2}/X^{0.49}.$$
Taking $X=(4000000\Npow2)^{50/49}$ and noting that $a_1=1$ with the
other terms on the left side being nonnegative, this says that
$e^{-1/10^6}\le RX^{1-b}\Gamma(1-b)+0.01$.
Since we assume $\Npow 2\ge 142$ we have $\log X\le 4.2\log\Npow2$,
and so we finally get that $X^{1-b}\le \exp(4.2/25)\le 1.19$.
We also have $\Gamma(1-b)\le 25\log\Npow2$,
and so we conclude that $R\ge {0.033\over \log\Npow2}$.
\end{proof}

\section{Isogenous curves, Manin constants,
and factors from twists and bad primes}
Finally we can turn to the other objects in our formula for the modular degree.
For the Manin constant we simply use the fact that $c\ge 1$.
So we have that
$${\rm deg}\,\phi\ge
{N\over\Omega}{0.033\over\log\Npow2}\cdot\prod_{p^2|N} U_p(1)^{-1}
\ge {ND^{1/6}\over 2675\log\Npow2}\cdot\prod_{p^2|N} U_p(1)^{-1}.$$

To bound the effects from the $U_p(1)^{-1}$, we recall its definition.
First we assume that $E$ is a global minimal twist; this is like requiring
that the model of $E$ be minimal at every prime, except that now we further
require that it be minimal when also considering quadratic twists.
See the author's paper \cite{watkins} for details.
For a minimal twist we have that $U_p(s)=(1-\epsilon_p/p^s)^{-1}$
where $\epsilon_p=+1,0,-1$, depending on certain properties
of inertia groups (see~\cite{watkins}).
In particular, we have that $\epsilon_p=+1$
for all primes congruent to 1~mod~12,
and $\epsilon_p=-1$ for all primes congruent to 11~mod~12.
When $p$ is 5~mod~12 we have that $\epsilon_p=+1$
exactly when $p^2|c_6$ and $p\|c_4$, while these conditions
imply that $\epsilon_p=-1$ for primes that are 7~mod~12.
Note in particular that $U_p(1)^{-1}$ is greater than~1
for primes that are 11~mod~12.
Also, when $U_p(1)^{-1}=1-1/p$ for primes that are 5~mod~12,
we have that $p^3|D$ while $p^2\|N$. For such primes we have that
$N_p(D_p)^{1/6}U_p(1)^{-1}\ge N_p^{7/6} p^{1/6}(1-1/p)\ge N_p^{7/6}$.
Finally we need to consider $p=2$ and~$p=3$.
There is not much to be done with $U_3(1)^{-1}$
except lower-bound it as $1-1/3=2/3$, whilst for $p=2$
in order for $U_2(s)$ to equal $(1-1/2^s)^{-1}$
we need that $2^8\|N$ and thus have $2^9\|D$.
So we have that $N_2(D_2)^{1/6}U_2(1)^{-1}\ge N_2^{7/6} 2^{1/6}(1-1/2)$.
Thus for a global minimal twist we have that
$${\rm \deg}\,\phi \ge {N^{7/6}\over 7150\log\Npow2}\cdot
\prod_{p^2|N\atop p\equiv 1\,(3)}(1-1/p).$$

We can estimate the product over primes using facts from prime number theory;
for $N\ge 20000$ the logarithm of the product is bounded by
$$\sum_{p^2|N\atop p\equiv 1\,(3)}{1\over p}+\sum_p {1/2p^2\over 1-1/p}
\le\sum_{p\le 1.02\log N\atop p\equiv 1\,(3)}{1\over p}+0.02
\le 0.5\log\log(1.02\log N)-0.33,$$
and so we have that
$${\rm \deg}\,\phi \ge {N^{7/6}\over 7150\log\Npow2}
{e^{0.33}\over\sqrt{0.02+\log\log N}}
\ge {N^{7/6}\over\log\Npow2}{1/5150\over\sqrt{0.02+\log\log N}}.$$

We wish to compare what happens on each side of this inequality upon
twisting our curve by an odd prime~$p$.
There are three cases, depending on the reduction type
of the minimal twist $F$ at~$p$.
If it has additive reduction, we simply have the same $U_p(s)$ as above,
and the conductor stays the same, with $D$ increasing by a factor of~$p^6$.
The modular degree goes up by~$p$ upon twisting, whilst the right side
of the inequality stays the same.
If $F$ has multiplicative reduction at~$p$,
we have that $U_p(s)^{-1}=(1-1/p^s)$. Here the conductor goes up by~$p$
and the discriminant by $p^6$ upon twisting, with the modular degree
gaining a factor of $(p^2-1)$. This is bigger than the factor of $p^{7/6}$
by which the right side increases.
Finally, if $F$ has good reduction at~$p$
we have that $U_p(2)^{-1}=(p-1)(p+1-a_p)(p+1+a_p)/p^3$ where $a_p$ is the
trace of Frobenius for~$F$. The conductor goes up by~$p^2$ and the
discriminant by~$p^6$.
The modular degree goes up by $(p-1)(p+1-a_p)(p+1+a_p)$, which is bigger
than the factor of $p^{7/3}$ gained by the right side,
even when $p=3$ and $a_p=\pm 3$.
So the above inequality is true for curves that are twist-minimal at~2.

Finally, we consider curves that are non-twist-minimal at~2.
If $2^8|N$, the right side of our inequality stays the same upon
twisting, while the left side does not diminish.
When 256 does not divide the conductor,
we can simply note that $U_2(1)^{-1}\ge (2-1)(2+1-2)(2+1+2)/2^3=5/8$
and so directly compute that
$${\rm deg}\,\phi\ge {ND^{1/6}\over 2675\log\Npow2}\cdot
{5\over 8}{2\over 3}\cdot{e^{0.33}\over\sqrt{0.02+\log\log N}}
\ge {N^{7/6}\over \log\Npow2}
{1/5000\over\sqrt{0.02+\log\log N}}.$$

\section{Summary of results}
We conclude this paper by giving a summary
of the various lower bounds that can be obtained.

\begin{theorem}
Suppose that $E$ is a rational semistable elliptic curve.
Then $${\rm deg}\,\phi_E\ge {N\over\Omega}{0.033\over 2\log N}
\ge {N^{7/6}\over 5350\log N}.$$
\end{theorem}

\begin{theorem}
Suppose that $E$ is a rational elliptic curve. Then
\begin{align*}
{\rm deg}\,\phi_E\ge {N\over\Omega}{0.033\over \log\Npow2}
\cdot\prod_{p^2|N} U_p(1)^{-1}&
\ge{N^{7/6}\over 7150\log\Npow2}\cdot
\prod_{p^2|N\atop p\equiv 1\,(3)} (1-1/p)\\
&\ge {N^{7/6}\over\log N}\cdot{1/10300\over\sqrt{0.02+\log\log N}}
\end{align*}
\end{theorem}

\begin{remark}
Note that we need $N\ge e^{86.8}$
in order to derive that ${\rm deg}\,\phi\ge N$ from Theorem~2.
\end{remark}

\end{document}